\documentclass[11pt]{amsart}
\usepackage{latexsym}
\usepackage{amssymb}
\usepackage{amsmath}
\usepackage{amsthm}

\newtheorem{theorem}{Theorem}[section]
\newtheorem{cor}[theorem]{Corollary}
\newtheorem{prop}[theorem]{Proposition}
\newtheorem{lemma}[theorem]{Lemma}

\theoremstyle{definition}

\theoremstyle{definition}

\theoremstyle{remark}

\theoremstyle{definition}

\newcommand{\comment}[1]{}
\newcommand{\R}{\mathbb{R}}
\newcommand{\T}{\mathbb{T}}
\newcommand{\Z}{\mathbb{Z}}
\newcommand{\Q}{\mathbb{Q}}
\newcommand{\coker}{\mathop{\textrm{coker}}}

\begin{document}

\title{3-Manifolds Admitting Toric Integrable Geodesic Flows}
\author{Christopher R. Lee}
\address{Department of Mathematics \\ University of Illinois \\ Urbana, IL 61801}
\email{crlee@math.uiuc.edu}
\date{\today}

\begin{abstract}
The geodesic flow of a Riemannian metric on a compact manifold $Q$ is said to be toric integrable if it is completely integrable and the first integrals of motion generate a homogeneous torus action on the punctured cotangent bundle $T^*Q\setminus{Q}$. If the geodesic flow is toric integrable, the cosphere bundle admits the structure of a contact toric manifold. By comparing the Betti numbers of contact toric manifolds and cosphere bundles, we are able to provide necessary conditions for the geodesic flow on a compact, connected 3-dimensional Riemannian manifold to be toric integrable. 
\end{abstract}

\maketitle

\section{Introduction} 
Let $(Q,g)$ be an $n$-dimensional Riemannian manifold. Consider, for each $t \in \R$, the diffeomorphism
\begin{eqnarray*}
\phi_{t}:TQ &\to& TQ \\
\phi_{t}(q,v) &=& (\gamma(t),\dot\gamma(t))
\end{eqnarray*}
where $\gamma$ is the unique geodesic with initial conditions $\gamma(0) = q$ and $\dot\gamma(0)=v$. Unique dependence of geodesics on initial conditions provides $\phi_{t} \circ \phi_{s} = \phi_{t+s}$, i.e. the family $\{\phi_{t}\}$ is actually a flow, the {\bf geodesic flow} of the metric $g$ on the manifold $Q$. 

Equivalently, we may describe the geodesic flow in terms of symplectic geometry. Recall that a {\bf symplectic manifold} is a manifold equipped with a closed, nondegenerate two-form $\omega$. The punctured cotangent bundle $T^*Q\setminus{Q}$ is a symplectic manifold; its symplectic form is $\omega = d\alpha$, where $\alpha$ is the one-form defined by $\alpha(v) = \eta(d\pi(v))$. Here, $\pi:T^*Q \to Q$ is the standard bundle projection, $\eta \in T^*Q$ and $v \in T_{\eta}(T^*Q)$. The form $\alpha$ is called the {\bf Liouville one-form}.

Define a smooth function $h:T^*Q\setminus{Q} \to \R$ by
\begin{displaymath}
h(q,p) = \sqrt{g^*_{q}(p,p)}
\end{displaymath} 
where $g^*$ is the metric on $T^*Q$ induced by $g$. The {\bf Hamiltonian vector field} $X_{h}$ of $h$ is then the unique solution to the equation
\begin{displaymath}
\iota_{X_{h}}\omega = dh.
\end{displaymath}
The geodesic flow of the metric $g$ on $TQ$ can then be identified with the flow of the Hamiltonian vector field $X_{h}$ on $T^*Q$. 

We say that the geodesic flow is {\bf toric integrable} if there exist $n$ smooth functions $h=f_1,\ldots,f_n$ on $T^*Q\setminus{Q}$ such that
\begin{enumerate}
\item $X_{f_i}(f_j)=0$ for all $i,j=1,\ldots,n$,
\item the set $\{df_i\}_{i=1}^{n}$ is linearly independent on an open, dense subset of $T^*Q\setminus{Q}$, 
\item each $f_i$ is homogeneous of degree 1: if $\rho_{t}(p,q) = (p,e^tq)$ is the $\R$-action given by dilations, then $\rho_{t}^*f_i = e^tf_i$ for each $i$, and
\item the $\R^n$-action defined by the flows of the vector fields $X_{f_{i}}$ descends to a well-defined action of $\T^n = \R^n/\Z^n$.  
\end{enumerate}
In the case that the geodesic flow is toric integrable, there is an effective action of the $n$-torus $\T^n$ on $T^*Q\setminus{Q}$ that preserves both the Liouville one-form $\alpha$ and the function $h$, and which also commutes with dilations.

Since the $\T^n$-action preserves $h$, it preserves to the {\bf cosphere bundle} $S(T^*Q)$ defined as the set of unit vectors in  $T^*Q\setminus{Q}$ with respect to the metric $g^*$. Recall that $S(T^*Q)$ is a {\bf contact manifold}, that is, the Liouville one-form $\alpha$ restricted to $S(T^*Q)$ satisfies $\alpha \wedge (d\alpha)^{n-1} \neq 0$. The $\T^n$-action generated by the integrals of the geodesic flow preserves $\alpha$, and so gives $S(T^*Q)$ the structure of a {\bf contact toric} manifold. The moral is thus:

\begin{center}
{\em If a Riemannian manifold $Q$ admits a toric integrable geodesic flow, then the cosphere bundle $S(T^*Q)$ of $Q$ is a contact toric manifold.}
\end{center}

The notion of toric integrability of the geodesic flow was important to Toth and Zelditch who discovered, while examining the spectra of Laplacians on compact Riemannian manifolds in \cite{tz}, that this condition provided substatntial simplification. However, toric integrability of the geodesic flow seems to be a rather restrictive geometric condition. It is known (see \cite{ls}) that the only metric which support a toric integrable geodesic flow on a torus is the flat metric. On the spheres $S^2$ and $S^3$, the geodesic flow of the round metrics are also toric integrable. Y. Colin de Verdiere in \cite{CdV} has constructed other metrics on $S^2$ which are toric integrable, but a complete classification of such metrics on low dimensional spheres is incomplete.  

We instead tackle a less formidable problem: for a given manifold $Q$, are there topological obstructions to toric integrability of the geodesic flow? In light of the preceeding discussion, this amounts to finding topological obstructions to cosphere bundles being contact toric manifolds. The main result of this paper is

\begin{theorem}
Suppose $Q$ is a $3$-dimensional compact, connected Riemannian manifold admitting a toric integrable geodesic flow. Then, the cosphere bundle $S(T^*Q)$ is diffeomorphic to either $S^3 \times S^2$ if $Q$ is simply connected or $\T^3 \times S^2$ if the fundamental group of $Q$ is infinite. If the fundamental group of $Q$ is finite and nontrivial, then $S(T^*Q)$ is homotopy equivalent to $(S^3/\Gamma) \times S^2$ where $\Gamma$ is a finite cyclic group. 
\end{theorem}

\section*{Acknowledgements} I would like to thank Eugene Lerman for his invaluable guidance and inspiration as well as Donald Yau for assistance with a calculation and beneficial conversations. In addition, Ian Hambleton provided some references and remarks concerning group actions on spheres for which I am grateful.

\section{Contact Toric Manifolds}
We begin by reviewing some basic facts about contact manifolds. Recall that a $(2n-1)$-dimensional contact manifold $M$ is {\bf cooriented} if there exists a globally defined contact form $\alpha$ on $M$. In this case, the {\bf contact structure} $\xi \subset TM$ is defined as the kernel of the contact form $\alpha$. We make the banket assumption in this paper that all contact manifolds in question are cooriented. If there is an effective action of $\T^n$ on $M$ preserving $\xi$, we say that $M$ is a {\bf contact toric} manifold. 

Recall that given a Hamiltonian action of a compact Lie group $G$ on a symplectic manifold $X$, we may associate to this action a (symplectic) moment map $\Phi:X \to \mathfrak{g}^{*}$ (here, $\mathfrak{g}^*$ is the dual of the Lie algebra of $G$). If $M$ is a contact toric manifold, consider the annihilator $\xi^{0} \subset T^*M$ of $\xi \subset TM$. The $\T^n$-action on $M$ lifts to a Hamiltonian $\T^n$-action on a component $\xi_{+}^{0}$. The {\bf contact moment map} is then defined to be the moment map $\Psi:\xi_{+}^{0} \to \mathfrak{t}^{*}$ of this action. Additionally, the image $\Psi(\xi_{+}^{0}) \cup \{0\}$ of the contact moment map together with the origin in $\mathfrak{t}^{*}$ is a rational polyhedral cone, called the {\bf moment cone} of the contact toric manifold $M$. 

The classification of compact contact toric manifolds was initiated by Banyaga and Molino in \cite{bm} and completed by Lerman in \cite{l2}. Lerman's result makes use of the moment cone of a contact toric manifold by singling out certain types of polyhedral cones in $\R^n$ he calls ``good cones'' and showing that they correspond to unique contact toric manifolds. We sum his results up in dimensions greater than 3:

\begin{theorem}{\cite{l2}}
Suppose $M^{2n-1}$ is a contact toric manifold with $n>2$. Then, 
\begin{enumerate}
\item if the $\T^n$-action is free, $M$ can be given the structure of a principal $\T^n$-bundle over $S^{n-1}$. Moreover, each principal $\T^n$-bundle over $S^{n-1}$ corresponds to a unique compact contact toric manifold. 
\item If the $\T^n$ action is not free, the moment cone of $M$ is good. Conversely, every good cone corresponds to a unique compact contact toric manifold.  
\end{enumerate}
\end{theorem}

We can say a bit more about the second part of the above classification theorem. If the moment cone of a contact toric manifold contains a subspace of positive dimension, then the contact toric manifold is diffeomorphic to $\T^{k} \times S^{2n-k-1}$ for some $k>0$. If the moment cone contains no lines, the contact toric manifold is said to be of {\bf Reeb type}.

Contact toric manifolds of Reeb type are not, in general, easily identifiable. To every contact manifold $(M,\alpha)$ we may associated a distinguished vector field, the {\bf Reeb vector field} $R_{\alpha}$ which is uniquely determined by the equations
\begin{eqnarray*}
\alpha(R_{\alpha}) &=& 1, \\
\iota_{R_{\alpha}}d\alpha&=&0.
\end{eqnarray*}
The Reeb vector field provides a splitting $TM = \xi \oplus \R{R_{\alpha}}$. Since the vector bundle $\xi$ is symplectic, there exists an almost complex structure $J$ which is compatible with the form $d\alpha|_{\xi}$ so that $\tilde{g} = d\alpha|_{\xi}(\cdot,J\cdot)$ is a metric on $\xi$. Since $TM$ splits, we can extend $\tilde{g}$ (by zero) to a metric on $TQ$. The Riemannian metric $g = \tilde{g} \oplus \alpha \otimes \alpha$ on $M$ is such that $\R{R_{\alpha}}$ and $\xi$ are orthogonal. If the Reeb vector field is Killing with respect to $g$, i.e. the Lie derivative of $g$ in the direction of $R_{\alpha}$ vanishes, we say that $M$ is a {\bf K-contact} manifold. What's important to us, as proved in Theorem 3.1 in \cite{l1}, is that a contact manifold is K-contact if and only if it is a contact toric manifold of Reeb type. 

Since we endeavor to study the geodesic flow of a compact, connected 3-dimensional Riemannian manifold, we should be interested in finding all compact contact toric manifolds in dimension 5. So far, our list includes
\begin{itemize}
\item principal $\T^3$-bundles over $S^2$,
\item $S^1 \times S^4$, $S^2 \times S^3$, and
\item K-contact manifolds.   
\end{itemize}
By the classification theorem, this list is exhaustive in dimension $5$. Luckily, thanks to Yamazaki, we also have at our disposal a complete list of all simply connected $5$-dimensional K-contact manifolds:

\begin{theorem}{\cite{yama}} 
Suppose $M$ is a closed, simply connected 5-manifold. Denote by $w_{2}(M)$ the second Stiefel-Whitney class of the bundle $TM \to M$. Then, $M$ admits an effective $\T^3$ action and a $\T^3$-invariant $K$-contact structure with $r$ closed orbits of the Reeb flow if and only if $M$ is diffeomorphic to
\begin{enumerate}
\item $S^5$, if $r=3$,
\item $\sharp^{r-3}(S^2 \times S^3)$, if $w_{2}(M)=0$, 
\item $(S^2 \tilde\times S^3) \sharp^{r-4} (S^2 \times S^3)$, if $w_2(M) \neq 0$. 
\end{enumerate}
where $S^2 \tilde\times S^3$ is the oriented, nontrivial $S^3$ bundle over $S^2$ and $\sharp$ denotes the connected sum operation.  
\end{theorem}

It is important to note here that if $M$ does indeed admit the structure of a contact toric manifold
with $r$ closed orbits of the Reeb flow, then $r \geq 3$. In particular, $r \geq 4$ if $w_{2}(M) \neq 0$.

\section{Cosphere Bundles and Simply Connected Contact Toric Manifolds}
In this section, we commence the proof of Theorem 1.1. We split the statement into cases, the first being the assumption that our Riemannian manifold $Q$ is simply connected. Of course, if $Q$ is simply connected, so is $S(T^*Q)$. This is easily seen by considering the long exact sequence of homotopy groups associated to the fibration $S^2 \to S(T^*Q) \to Q$. In particular, we prove
\begin{prop}
Suppose $(Q,g)$ is a $3$-dimensional compact, simply connected Riemannian manifold admitting a toric integrable geodesic flow. Then, the cosphere bundle of $Q$ is diffeomorphic to $S^3 \times S^2$.  
\end{prop}

The crux of our computations involve the following theorem, due to McCord, Meyer and Offin (\cite{mmo}) which place necessary conditions on a manifold to be a cosphere bundle. We adopt the following notation: $\beta_{i}(M)$ is the $i$th Betti number of $M$ and $T_{i}(M)$ is the torsion subgroup of the homology group $H_{i}(M,\Z)$. 
\begin{theorem}\cite{mmo}
Suppose $M$ is a connected, $(2n-1)$-dimensional manifold. Then, if $M$ arises as the cosphere bundle of some $n$-manifold $Q$, 
\begin{itemize}
\item[(i)] $\beta_{2n-1}(M)=1$ and $T_{n-1}(M)$ is cyclic or trivial if $Q$ is compact and orientable,
\item[(ii)] $\beta_{2n-1}(M)=1$ and $T_{n-1}(M)$ has order $4$ if $Q$ is compact and nonorientable.
\end{itemize}
In addition, if $T_{n-1}(M)$ is trivial and $M$ is the cosphere bundle of some compact, orientable manifold $Q$, we have
\begin{eqnarray}
0 \leq \beta_{n-1}(M)-\beta_{2n-2}(M) &=& \beta_{n}(M)-\beta_{1}(M) \leq 1 \\
\left|\sum_{i=0}^{n-2} (-1)^{i}\beta_{i}(M)+(-1)^{n}(1-\beta_{2n-2})\right| &=& 1+\beta_{1}(M)-\beta_{n}(M). 
\end{eqnarray}
\end{theorem}

From the results of the previous section, namely Theorem 2.2, our candidates for simply connected, contact toric cosphere bundles are $S^5$and connected sums of $S^2 \times S^3$ and $S^3 \tilde\times S^2$. Our first application of Theorem 3.2 will be to show that $S^5$ is not the cosphere bundle of any manifold $Q$. We then explicitly compute the Betti numbers of the remaining simply connected contact toric manifolds on our list. 
\begin{lemma}
The contact toric manifold $S^5$ does not arise as the cosphere bundle of any compact, simply connected 3-manifold.
\end{lemma}
\begin{proof}
Suppose that $S^5$ arose as the cosphere bundle of a compact manifold. Then, by Theorem 3.2, we would have 
\begin{displaymath}
|\beta_{0}(S^5)-\beta_{1}(S^5)-1+\beta_{4}(S^5)| = 1+\beta_{1}(S^5)-\beta_{3}(S^5).
\end{displaymath}
But, $\beta_{k}(S^5)$ is trivial for $0<k<5$, and $\beta_{0}(S^5)=1$. Hence, we would then have $0=1$, a contradiction. 
\end{proof}

\begin{lemma}
For $k>0$, we have
\item 
\begin{displaymath}
\beta_{i}(\sharp^{k}(S^2 \times S^3))= \left\{ \begin{matrix}  1 & \textrm{if } i=0,5 \\ k & \textrm{if } i=2,3 \\ 0 & \textrm {if } i= 1,4 \textrm{ or } i >5  \end{matrix} \right.
\end{displaymath}
\end{lemma}
\begin{proof}
The case $k=1$ follows from applying the K\"unneth formula to $S^2 \times S^3$ to get Betti numbers
\begin{displaymath}
\beta_{i}(S^2 \times S^3) = \left\{ \begin{matrix} 1 & \textrm{if } i=0,2,3,5 \\ 0 & \textrm {if } i= 1,4 \textrm{ or } i >5 \end{matrix} \right.
\end{displaymath}

For $0<i<m$ and connected, orientable $m$-manifolds $M,N$, we have $H_{i}(M \sharp N) \cong H_{i}(M) \oplus H_{i}(N)$. Also, $\dim_{\Z}{H_0(M \sharp N)} = 1$ by connectivity and $\dim_{\Z}{H_m(M \sharp N)} = 1$ by orientability. Our result then follows by induction. 
\end{proof}

\begin{lemma}
For $k \geq 0$, 
$$
\beta_{i}((S^2 \tilde\times S^3) \sharp (\sharp^{k}(S^2 \times S^3)))= \left\{ \begin{matrix}  1 & \textrm{if } i=0,5 \\ k+1 & \textrm{if } i=2,3 \\ 0 & \textrm {if } i= 1,4 \textrm{ or } i >5  \end{matrix} \right.
$$
\end{lemma}
\begin{proof}
When $k=0$, we are making a statement about the Betti numbers of $S^2 \tilde\times S^3$. Now, there is an exact sequence in homology, called the Wang sequence, which is useful for fibrations $F \hookrightarrow E \to S$ whose base is a sphere: 
$$
\cdots\to H_{m}(F) \to H_{m}(E)\to H_{m-n}(F) \to H_{m-1}(F) \to \cdots.
$$
The Wang sequence for the fiber bundle $S^3 \hookrightarrow S^2 \tilde\times S^3 \to S^2$ thus yields
$$
\beta_{i}(S^2 \tilde\times S^3) = \left\{ \begin{matrix} 1 & \textrm{if } i=0,2,3,5 \\ 0 & \textrm {if } i= 1,4 \textrm{ or } i >5 \end{matrix} \right.
$$
Induction and Lemma 3.4 now combine to complete the argument. 
\end{proof}

The previous results narrow down possible candidates for contact toric cosphere bundles:

\begin{cor}
The only possible simply connected, compact contact toric cosphere bundles are $S^2 \times S^3$ and $S^2 \tilde\times S^3$. 
\end{cor}

\begin{proof}
When $w_2(M) =0$ and $r \geq 3$, equation $(2)$ of Theorem 3.2 applied to $\sharp^{r-3} (S^2 \times S^3)$ gives $0=4-r$.
By equation $(1)$ of the same theorem, the only possibilities for $r$ are $r=3$ or $r=4$. But $0 \neq 4-3= 1$, and
when $r=4$ we have $S^2 \times S^3$.

Similarly, when $w_2(M) \neq 0$, we have $r \geq 4$. Equation $(1)$ again provides $0=4-r$. This
is only possible when $r=4$ and hence $S^2 \tilde\times S^3$ is a possibility.  
\end{proof}

We are now in a position to prove Proposition 3.1. Before we do this, however, note that the Lie group $S^3$ is parallelizable. Indeed the cosphere bundle is diffeomorphic to $S^3 \times S^2$. To complete the proof of Proposition 3.1, then, it suffices to show that $S^3 \tilde\times S^2$ cannot arise as the cosphere bundle of any manifold. 

\begin{proof}{(of Proposition 3.1)}
Suppose $S^3 \tilde\times S^2$ is the cosphere bundle of some compact Riemannian manifold $Q$. Note that, by Theorem 3.2(i), $Q$ must be orientable since $T_2(S^2 \tilde\times S^3)$ is trivial. So, $Q$ is compact and orientable and by a theorem of Stiefel (see \cite{parker} for a proof using obstruction theory) it is parallelizable. Hence, $TQ = Q \times \R^3$. So, $S(TQ) = Q \times S^2$.  Again, since $Q$ is parallelizable, $w_2(Q \times S^2) = w_2(S^2)$. But, $w_2(S^2)=0$ by naturality of the Stiefel-Whitney classes (see \cite{milnor} for details). Since $w_{2}(S^2 \tilde\times S^3) \neq 0$, it cannot occur as the cosphere bundle of $Q$. 
\end{proof}

\section{Non-Simply Connected Contact Toric Cosphere Bundles}
In the previous section, we found that the only compact, simply connected contact toric manifold which is also a cosphere bundle is $S^3 \times S^2$. We now turn our attention to the non-simply connected case. 

Recall that $\T^n$ with the flat metric always has a toric integrable geodesic flow. As it is a Lie group, $\T^n$ is parallelizable and its cosphere bundle is $\T^n \times S^{n-1}$. From Theorem 2.1, we see that $\T^n \times S^{n-1}$ is a contact toric manifold. The main result of this section is

\begin{prop}
Suppose $Q\neq \T^3$ is a compact, connected Riemannian manifold with nontrivial, finite fundamental group. Then, if $Q$ admits a toric integrable geodesic flow, the cosphere bundle of $Q$ is homotopy equivalent to the cosphere bundle of $S^3/\Gamma$, where $\Gamma$ is a finite cyclic group. 
\end{prop}

The proof of Proposition 4.1 will proceed as follows: We first show that if $Q \neq \T^3$ and the cosphere bundle $S(T^*Q)$ is contact toric, then the fundamental group of $S(T^*Q)$ is finite. We then show that if $M$ is a compact, connected contact toric manifold which arises as a cosphere bundle, the universal cover $\widetilde{M}$ is also a compact contact toric manifold which is a cosphere bundle. Since universal covers are simply connected, we apply the results of the previous section.

\begin{lemma} Suppose that $M$ is a non-simply connected 5-manifold, $M = S(T^*Q)$ for some compact, connected
3-manifold $Q \neq \T^3$. Then, the fundamental group of $M$ is finite. 
\end{lemma} 

\begin{proof}
In Theorem 2.2, we see that the only $5$-dimensional contact toric manifolds which are not simply connected and have infinite fundamental group are diffeomorphic to principal $\T^3$ bundles over $S^2$, $\T^2 \times S^3$ or $S^1 \times S^4$. Now, applying the K\"unneth formula to the latter two spaces, we see that $\beta_{2}(\T^2 \times S^3)-\beta_{4}(\T^2 \times S^3) = 1-2 <0$ and $\beta_{2}(S^1 \times S^4)-\beta_{4}(S^1 \times S^4) = 0-1 <0$. So, by Theorem 3.2, neither of these spaces arise as the cosphere bundle of any manifold. 

It remains to show that a (nontrivial) principal $\T^3$-bundle $P$ over $S^2$ cannot arise as the cosphere bundle over any compact 3-manifold $Q$. In order to calculate the Betti numbers of such bundles, we employ the cohomology Serre spectral sequence with rational coefficients. Denote by $E_{n}^{p,q}$ the $(p,q)$th entry on the $n$th page of the spectral sequence. Denote by $a,b,c$ the generators of $H^1(\T^3;\Q)$ and by $x$ the generator of $H^2(S^2;\Q)$. Then, using the multiplicative structure on the spectral sequence we have
$$
\begin{array}{llll}
E_{2}^{0,0} = \Q & E_{2}^{0,1} = \Q{a}\oplus\Q{b}\oplus\Q{c} & E_{2}^{0,2} = \Q{ab}\oplus\Q{ac}\oplus\Q{bc} & E_{2}^{0,3} = \Q{abc} \\
E_{2}^{2,0} = \Q{x} & E_{2}^{2,1} = \Q{ax}\oplus\Q{bx}\oplus\Q{cx} & E_{2}^{2,2} = \Q{abx}\oplus\Q{acx}\oplus\Q{bcx} & E_{2}^{2,3} = \Q{abcx}
\end{array} 
$$
and all other entries of the sequence are 0. Thus, the only nonzero differentials are
\begin{eqnarray*}
d_{2}^{1}:E_{2}^{0,1} \to E_{2}^{2,0} \\
d_{2}^{2}:E_{2}^{0,2} \to E_{2}^{2,1} \\
d_{2}^{3}:E_{2}^{0,3} \to E_{2}^{2,2}
\end{eqnarray*}

For the generators $a,b,c$ of $H^1(\T^3;\Q)$, we write $d_{2}^{1}(a) = a'x$, $d_{2}^{1}(b) = b'x$, $d_{2}^{1}(c) = c'x$ where
$a',b',c' \in \Q$. We then have
\begin{eqnarray*}
d_{2}^{2}(ab) &=& d_{2}^{2}(a)b-ad_{2}^{2}(b) \\
&=& a'bx-b'ax 
\end{eqnarray*}
and similarly, $d_{2}^{2}(ac) = a'cx-c'ax$, $d_{2}^{2}(bc) = b'cx-c'bx$. In addition,
\begin{eqnarray*}
d_{2}^{3}(abc) &=& d_{2}^{3}(a)bc-ad_{2}^{3}(bc) \\
&=& a'bcx-b'acx+c'abx
\end{eqnarray*}

Now, we may break our computations down into essentially 4 cases: all of $a',b',c'$ are zero, none are zero, at most one is zero and at most two are zero. These cases cover all possible situations as the above expressions represent a symmetric situation. Now, when $a'=b'=c'=0$, the principal bundle is trivial so we may assume at least one is nonzero.

For illustration, we now compute the Betti numbers of $P$ when none of $a',b',c'$ are zero. In this
case, $E_{\infty}^{p,q} = E_{3}^{p,q}$, $E_{3}^{0,i} = \ker{d_{2}^{i}}$ and $E_{3}^{2,i} = \coker{d_{2}^{i}}$.
Since we are only looking for Betti numbers, we need only find the sums of the ranks of kernels and 
cokernels along the ``diagonals''. Now, $\ker{d_{2}^{1}}$ contains terms like $b'c-c'b$ and is hence free
on 3 generators. So, $\beta_{1}(P) = 3$ since $E_{3}^{1,0} =0$. Similarly, we can see that $\beta_2(P) = 1$, 
$\beta_{3}(P) = 2$, $\beta_{4}(P) = 2$ and $\beta_{5}(P) = 1$. By Theorem 3.1 (2), we see that in this case 
$P$ cannot arise as the cosphere bundle of a compact, connected manifold $Q$. 

In a similar fashion, we may compute the Betti numbers of $P$ in each of the above cases. As it turns out, 
when exactly one or two of $a',b',c'$ is zero we have
$$
\beta_{i}(P) = \left\{ \begin{array}{cc} 2 & \textrm{if } i=1,4\\  1 & \textrm{if } i=0,2,3,5\end{array}\right.
$$
which again violates conditions imposed in Theorem 3.2. It follows that none of the possible total
spaces $P$ under consideration arise as cosphere bundles.     
\end{proof}

\begin{lemma}
Suppose $M = S(T^*Q)$ for some manifold $Q$. Then, if $\widetilde{M}$ is the universal cover of $M$, $\widetilde{M} = S(T^*\widetilde{Q})$
where $\widetilde{Q}$ is the universal cover of $Q$.  
\end{lemma}
\begin{proof}
The covering map $p:\widetilde{Q} \to Q$ is a local diffeomorphism. Hence, it induces an isomorphism on fibers $S(T^*_{q}Q) \to S(T^*_{p(q)}\widetilde{Q})$. Since $\widetilde{Q}$ is simply connected, so is $S(T^*\widetilde{Q})$. It follows that $\widetilde{M} = S(T^*\widetilde{Q})$. 
\end{proof}

\begin{lemma}
Suppose $M$ is a contact toric manifold with finite fundamental group. Denote by $\widetilde{M}$ the universal cover of $M$. Then, $\widetilde{M}$ is a compact contact toric manifold. 
\end{lemma} 
\begin{proof}
We have to verify both that $\widetilde{M}$ admits a contact structure and that the torus action on $M$ lifts to a torus
action on $\widetilde{M}$. 
First note that being a contact manifold is a local condition. Since a covering map is a local diffeomorphism, 
ths contact structure on $M$ lifts to a contact structure on $\widetilde{M}$. 

We now lift the torus action from $M$ to $\widetilde{M}$. We claim it is enough to lift an $S^1$ action
since a product of groups acts on a manifold if and only if each component does and these actions commute. 

So, suppose $S^1$ acts smoothly on $M$. This action generates a vector field $X$ on $M$. Denote by $\phi_{t}$ the
flow of $X$. If we view $S^1$ as $\R/\Z$, we see that $\phi_{t} = \phi_{t+1}$ for any $t$. Since
the covering map $\pi: \widetilde{M} \to M$ is a local diffeomorphism, we may lift $X$ to a vector field $\widetilde{X}$ 
on $\widetilde{M}$ with corresponding flow $\tilde{\phi}$. We claim that for some positive integer
$k$, $\phi_{t} = \phi_{t+k}$ for all $t$. Hence, the $S^1$ action will lift to $\widetilde{M}$ (with perhaps
a different parameterization of $S^1$). Now, the fundamental group of $M$ is finite, so for any point $x \in M$, the fiber
$\pi^{-1}(x)$ is finite, say cardinality $m$. Since $M$ is assumed connected, the number of elements
in the fiber is constant over each point in $M$ so we need only argue pointwise. Now, for any $\tilde{x} \in \pi^{-1}(x)$, 
note that $\tilde{\phi}(\tilde{x}) \in \pi^{-1}(x)$ as well since $\tilde{\phi}$ was gotten by lifting
$\phi$ locally. In other words, $\pi(\tilde{\phi}(\tilde{x})) = \phi(x)$. Since the fibers have only finitely
many elements, and the flow of $\tilde{X}$ always stays within fibers, there must be some $k$
with $\tilde{\phi}_{t}(\tilde{x}) = \tilde{\phi}_{t+k}(\tilde{x})$ for all $t$.  
\end{proof}

We now provide a proof of Proposition 4.1. which, in turn, completes the proof of Theorem 1.1. 
\begin{proof}{(Proposition 4.1)}
Let $M = S(T^*Q)$. Then, by Lemmas 4.3 and 4.4, the universal cover $\widetilde{M}$ is a contact toric manifold and the cosphere bundle of $\widetilde{Q}$. Since $\widetilde{M}$ is simply connected, Proposition 3.1 implies that $\widetilde{M}$ is diffeomorphic to $S^3 \times S^2$. Since $\widetilde{Q}$ is a closed $3$-manifold that is simply connected, Whitehead's theorem (along with the Hurewicz theorem) implies that $\widetilde{Q}$ is homotopy equivalent to $S^3$. Therefore, $Q$ is homotopy equivalent to $S^3/\pi_1(Q)$. Now, consider the bundle $S^2 \hookrightarrow M \to Q$. Since $S^2$ is simply connected, by the long exact sequence in homotopy we can see that $\pi_1(M) \cong \pi_1(Q)$. Hence, $Q$ is homotopy equivalent to $S^3/\pi_1(M)$. By passing to a double cover, if necessary, we may assume without loss of generality that $Q$ is orientable. Hence, it is parallelizable and the cosphere bundle of $Q$ is homotopy equivalent to $S^3/\Gamma \times S^2$. 

In \cite{l1} it is shown that the fundamental group of a contact toric $\T^n$-manifold of Reeb type is a quotient of $\Z^n$ by a syblattice. So,  $\pi_1(M)$ is isomorphic to a product of at most 3 finite cyclic groups. 

We now show that $\Gamma = \pi_1(M)$ is indeed a cyclic group. Write $\Gamma=\Z_{m} \times \Z_{l} \times \Z_{k}$. If the greatest common divisor $(m,l,k) :=\textrm{gcd}(m,l,k)=1$, then $\Gamma$ is isomorphic to the finite cyclic group $\Z_{(m,l,k)}$. So, suppose not. Then, $\textrm{gcd}(m,l,k)=s \neq 1$ and there is some prime $p$ dividing $s$ and hence each of $m,l,k$. So, $\Gamma$ contains $\Z_p \times \Z_p \times \{0\}$ as a noncyclic subgroup of order $p^2$. A theorem of Madsen,Thomas and Wall \cite{mtw} says that a finite group cannot act freely on $S^3$ if it contains a noncyclic subgroup of order $p^2$. Hence,  $\Gamma$ cannot act freely on $S^3$ and cannot be an isomorphic copy of $\pi_1{M}$.  A similar argument applies if $\Gamma$ is assumed to be  a product of two finite cyclics, and so the only way $\Gamma$ can act freely on $S^3$ is if it is a finite cyclic group. 
\end{proof}

\end{document}